\documentclass[oneside,11pt]{article}
\usepackage{amsfonts}
\usepackage{amsmath}
\usepackage{graphicx}

\begin{document}

\title{Regular level sets of Lyapunov graphs of nonsingular Smale flows on 3-manifolds \textit{} }

\author {Bin Yu }

\maketitle

\hspace*{-0.2cm} \textbf{Abstract:}
 In this paper, we first discuss the regular level set of a nonsingular
 Smale flow (NSF) on a 3-manifold. The main result about this topic is
 that a 3-manifold $M$ admits an NSF flow which has a regular level set homeomorphic to
 $(n+1)T^{2}$  $(n\in \mathbb{Z}, n\geq 0)$
      if and only if $M=M'\sharp n S^{1}\times S^{2}$.  Then we
      discuss how to realize a template as a basic set of an NSF on
      a 3-manifold. We focus on the connection between the genus of the
      template $T$  and the topological structure of the realizing
      3-manifold $M$.

\hspace*{-0.2cm} \textbf{Keywords:} Lyapunov graphs, nonsingular
Smale flows, three manifolds, templates.

\hspace*{-0.2cm} \textbf{MR(2000)Subject Classification:} 57N10,
58K05, 37E99, 37D45.

 \vspace*{0.5cm}
\begin{bfseries}
1.Introduction
\end{bfseries}
\vspace*{0.5cm}

  The qualitative study of a class of smooth flows on a
       manifold is often done by understanding the behavior of
       the flow on the basic sets and describing the disposition of
       the basic sets on the manifold.

       Smale flow is an important class of flow. There are many
       works about Smale flows on 3-manifolds. F.B$\acute{e}$guin and
       C.Bonatti [BB] described the behavior of the neighborhood of
       a basic set of an Smale flow on a 3-manifold. M.Sullivan [Su]
       used template and knot theory to describe more embedding
       information of a special type of nonsingular Smale flows on
       $S^{3}$. J.Franks ([F3], [F4]) used homology to describe some
       embedding information of nonsingular Smale flows (NSF) on
       3-manifolds.

      In particular, J.Franks [F1] introduced Lyapunov graph to give
     a global picture of how the basic sets of an NSF on $S^{3}$ are
     situated on $S^{3}$. Following the idea of J.Franks, K. de Rezende [R]
     used Lyapunov graph to classify Smale flows on $S^{3}$. For a
     smooth flow $\phi_{t}:M\rightarrow M$ with a Lyapunov function
     $f:M\rightarrow R$, the Lyapunov graph $L$ is a rather natural
object. The idea is to construct an oriented graph by identifying to
a point each component of $f^{-1}(c)$ for each $c\in R$. For Smale
flow, the weight of a edge of the Lyapunov graph is defined to be
the genus of the regular set of the edge, see [R].

 It is interesting to generalize J.Franks' work to all
 3-manifolds. This amounts to determine the necessary and sufficient condition
 when an abstract L graph is associated with an NSF on a
 3-manifold. N.Oka made an attempt on this topic in his paper
 [Ok].

 To generalize J.Franks' work, we should study how the topology of
$M$ is related to the Lyapunov graph $L$ of a Smale flow on $M$. In
particular, we are interested in studying: (1) the relation between
the topology of $M$ and the topology of $L$; (2) the relation
between the topology of $M$ and the weights of edges of $L$.

For $M=S^{3}$, J.Franks proved that the Lyapunov graph of any NSF on
$S^{3}$ must be a tree and the weight of any edge of $L$ is 1, i.e.,
 the regular level sets of any NSF on $S^{3}$ must be tori. See [F1]. For general manifold $M$,  R.N.Cruz and K. de Rezende [CR]
studied the relation between the topology of $M$ and the topology of
$L$. In particular, they showed that the cycle-rank
   $r(L)$ of the Lyapunov graph $L$ of any Smale flow on a 3-manifold $M$ satisfies: $r(L)\leq g(M)$.
   Here $g(M)$ is the maximal number of mutually disjoint, smooth, compact,
   connected, two-sided codimension one sub-manifolds that do not disconnect $M$.

 In the first part of this paper, we will discuss the connection between the weights of edges of
 $L$ and the topological structure of the 3-manifold. By
   [F1] and
   [R], for any integer $n\geq 0$ and any 3-manifold $M$, there exists a Smale flow
   such that the genus of a regular level set
   of a Lyapunov graph associated with the flow is $n$. Therefore in some sense the
   interesting part is NSF. The following theorem is related to this
   topic.\\

      \textbf{Theorem 1}
      \begin{itshape}
      Suppose a closed orientable 3-manifold M admits
      a nonsingular Smale flow $\varphi_{t}$ with Lyapunov function
      $f$. Let $L$ be the Lyapunov graph of $f$. Then M admits an NSF flow which
      has a regular level set homeomorphic to
 $(n+1)T^{2}$  $(n\in \mathbb{Z}, n\geq 0)$
      if and only if $M=M'\sharp n S^{1}\times S^{2}$. Here
      $(n+1)T^{2}$ means the connected sum of $n+1$ tori and
      $M'$ is any closed orientable 3-manifold. The cyclic number of
      a Lyapunov graph of the NSF, $r(L)\geq n$.
      \end{itshape}\\

      Another interesting problem is to discuss how to realize
      1 dimensional basic sets in an NSF on a 3-manifold. This will be the topic of the second part of this paper.

      If the basic
      sets are described by a suspension of a subshift of finite
      type (SSFT), the realization problem is considered in [PS] and
      [F2]. In particular, by the result of [F1] and Proposition 6.1 of [F4], it is easy to show that any SSFT can be
      realized as a basic set of an NSF of any closed orientable
      3-manifold.

      If the basic sets are described by template, Theorem 3.5.1 in Meleshuk [Me]
      described an algorithm to decide whether an embedded template can
      be  realized as a basic set of an NSF on $S^{3}$. Meleshuk [Me] also proved (see [Me] Theorem 3.4.1)
      that for any template there exists some 3-manifold with an NSF on it having a basic set modeled by the template.
      M.Sullivan [Su] and the author [Yu] considered how to realize an NSF with three
    basic sets: a repelling orbit $r$, an attracting orbit $a$, and a
     non-trivial saddle set modeled by a Lorenz template and a Lorenz
     like template on 3-manifolds. G.Frank [Fr] gave some obstruction to the realization of certain templates
     in a homology 3-sphere. The first statement of Theorem 2 is a generalization of Theorem A in [Fr].
     The second statement is a direct corollary of the first one and
     Theorem 3.4.1 of [Me]. But we will give a constructive proof of the second statement in section 4. Actually our construction
      also proves that for any template there exists some 3-manifold with an NSF on it having a basic set modeled by the template. \\

      \textbf{Theorem 2}
      \begin{itshape}
      Let $T$ be a template .

     (1).  If a closed orientable 3-manifold $M$ admits an NSF with a basic set modeled
     by $T$, then $M=M'\sharp g(T) S^{1}\times S^{2}$.  Moreover,  $r(L)\geq g(T)$,
      where $r(L)$ is the cyclic number of the Lyapunov graph $L$ of the NSF.

     (2). There exists a closed orientable
     3-manifold $M'$  such that
     $M=M'\sharp g(T) S^{1}\times S^{2}$ admits an NSF with a basic set
     modeled by $T$.
     \end{itshape}\\

     In Theorem 2, $g(T)$ is a number which describes the basic set modeled by $T$,
     see Definition 2.7. In Theorem 1 and Theorem 2 above, $S^{1}\times S^{2}$ plays an important role. The similar phenomenon appears
 in the study of nonsingular Morse-Smale flows on 3-manifolds, see  M. Saito's paper [Sa].

     The paper is organized as follows. In Section 2, we give some definitions and detailed background
     knowledge. Theorem 1 and Theorem 2 are proved in Section 3 and
     Section 4 respectively. In Section 5, we use thickened template
     and some surgeries to give a visualization of an NSF with $2T^{2}$ regular
     level set on $S^{1} \times S^{2}$. In Section 6, we discuss possible development
     of this work by asking two questions.

\vspace*{0.5cm}
\begin{bfseries}
2. Preliminaries
\end{bfseries}
\vspace*{0.5cm}

For the standard material of hyperbolic dynamical systems, we refer
the reader to the book written by C. Robinson [Ro]. For a more
detailed
account of Smale flows, see [F3].\\

\textbf{Definition 2.1} A smooth flow $\phi_{t}: M\rightarrow M$ on
a compact manifold is called a \emph{Smale flow} if:

     (1) the chain recurrent set $R$ has hyperbolic structure;

     (2) the $dim(R)\leq 1$;

     (3) $\phi_{t}$ satisfies the transverse condition.

  If a Smale flow $\phi_{t}$ has no singular point, we call
  $\phi_{t}$ \emph{nonsingular Smale flow} (NSF).\\

\textbf{Definition 2.2} A Lyapunov graph $L$ for a flow $\phi_{t}:
M\rightarrow M$ and a Lyapunov function $f:M\rightarrow R$ is
obtained by taking the quotient complex of $M$ by identifying to a
point each component of a level set of $f$. Denote the cyclic
number (or the first Betti number) of $L$  by $r(L)$.\\

\textbf{Definition 2.3} An \emph{abstract Lyapunov graph} is a
finite, connected and oriented graph $L$ satisfying the following
two conditions:

(1) $L$ possesses no oriented cycles;

(2) each vertex of $L$ is labeled with a chain recurrent flow on a
compact space. \\

\textbf{Theorem 2.4 (Bowen [Bo])}
\begin{itshape}
If $\phi_{t}$ is a flow with hyperbolic chain recurrent set and
$\Lambda$ is a 1-dimensional basic set, then $\phi_{t}$ restricted
to $\Lambda$ is topologically equivalent to the suspension of a
basic subshift of finite type (i.e., a subshift
associated to an irreducible matrix).
\end{itshape}\\

This theorem enables us to label a vertex of the (abstract) Lyapunov
graph that represent the 1-dimensional basic sets of a Smale flow
with the suspension of an SSFT $\sigma(A)$. For simplicity we will
label the vertex with the nonnegative integer irreducible matrix
$A$. For a vertex labeled with matrix $A=(a_{ij})$, let
$B=(b_{ij})$, where $b_{ij}\equiv a_{ij} (mod 2)$ and $k=\dim
\ker((I-B):F_{2}^{m}\rightarrow F_{2}^{m})$, $F_{2}=\mathbb{Z}/2$.
The number of incoming (outgoing) edges is denoted by $e^{+}$
($e^{-}$). Denote the weight of an edge of the Lyapunov graph by the
genus of the regular level set of the edge. Let
$g_{j}^{+}(g_{j}^{-})$ be the weight on an incoming (outgoing) edge
of the vertex. In this paper, we always denote the initial vertex
and the terminal vertex of an oriented edge $E$ by $i(E)$ and $t(E)$
respectively. Whenever we talk about the Lyapunov graph $L$ of an
NSF on a 3-manifold $M$, L is always associated with a Lyapunov
function $f:M\rightarrow R$ and a map $h:M\rightarrow L$ such that
$f=\pi \circ h$. Here $\pi: L\rightarrow R$ is the natural
projection. By the work of F.B$\acute{e}$guin and C.Bonatti [BB], it
is easy to show that there is a unique Lyapunov graph $L$ of a Smale
flow on a 3-manifold.

The following classification theorem is due to K. de Rezende [R]:\\

\textbf{Theorem 2.5}
\begin{itshape}
Given an abstruct Lyapunov graph $L$ whose sink (source) vertices
are each labeled with an attracting (repelling) periodic orbit or an
index 0 (index 3) singularity, then $L$ is associated with a Smale
flow $\phi_{t}$ and a Lyapunov function $f$ on $S^{3}$ if and only
if the following conditions hold:

(1) The underlying graph $L$ is a tree with exactly one edge
attached to each sink or source vertex.

(2) If a vertex is labeled with an SSFT with matrix $A_{m\times m}$,
then we have
$$
\begin{array}{c}
e^{+}> 0, e^{-}> 0,\\
k+1-G^{-}\leq e^{+}\leq k+1, \text{\,with\,}
G^{-}=\sum_{i=1}^{e^{-}}g_{i}^{-} and\\
k+1-G^{+}\leq e^{-}\leq k+1, \text{\,with\,}
G^{+}=\sum_{j=1}^{e^{+}}g_{j}^{+}.
\end{array}$$

(3) All vertices must satisfy the Poincare-Hopf condition. Namely,
if a vertex is labeled with a singularity of index r, then

$$(-1)^{r}=e^{+}-e^{-}-\sum g_{j}^{+}+ \sum g_{i}^{-}$$

and if a vertex is labeled with a suspension of an SSFT or a
periodic orbit, then

$$0=e^{+}-e^{-}-\sum g_{j}^{+}+ \sum g_{i}^{-}$$.
\end{itshape}\\

Template theory was first introduced to dynamics by R.F.Williams and
J.Birman in their papers [BW1], [BW2]. [GHS] is a monograph on this
subject. As a model of a basic set of an NSF, it provides more information than SSFT. \\

\textbf{Definition 2.6} A $Template$  $(T,\phi)$ is a smooth
branched 2-manifold $T$, constructed from two types of charts,
called $joining$ $charts$ and $splitting$ $charts$, together with a
semi-flow. A semi-flow is the same as flow except that one cannot
back up uniquely. In Figure 1 the semi-flows are indicated by arrows
on charts. The gluing maps between charts must respect the semi-flow
and act linearly on the edges.

\begin{center}

\includegraphics[totalheight=5cm]{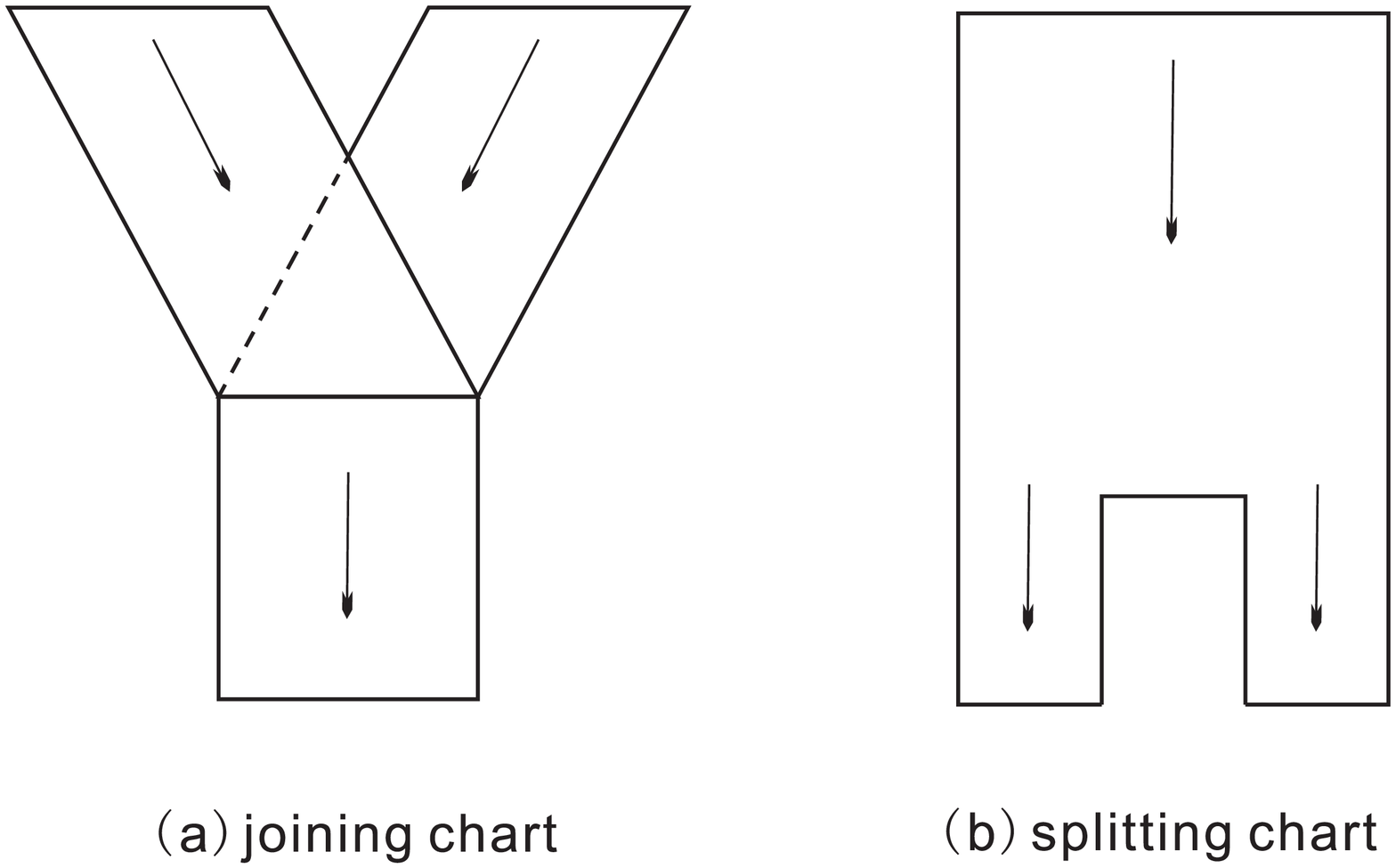}

\begin{center}
Figure 1
\end{center}
\end{center}

If we extend a template in the direction perpendicular to its
surface, we have a thickened template. The inverse limit of the
semi-flow on the template produces
 a flow on the thickened template (see [Me] or [Su]), as Figure 2 shows. For a
template $T$, there exists a unique \emph{thickened template}
$\overline{T}$. $\partial \overline{T}$ is composed of entrance sets
$X$, exit sets $Y$ and dividing curves. Let $X$ be a union of $s$
components: $X_{1}, ..., X_{s}$ and $Y$ be a union of $t$
components: $Y_{1}, ..., Y_{t}$. Suppose the genus of $X_{i}$ is
$n_{i}^{+}$ and the genus of $Y_{j}$ is $n_{j}^{-}$. Reindexing if
necessary, we may assume

 $\left\{
  \begin{array}{ll}
    n_{i}^{+}> 1, & \hbox{$i=1,...,s_{0}$;} \\
    n_{i}^{+}=0, & \hbox{$i=s_{0}+1,...,s_{0}+s_{1}$;} \\
    n_{i}^{+}=1, & \hbox{$i= s_{0}+s_{1}+1,...,s$.}
  \end{array}
\right.$

and

$\left\{
  \begin{array}{ll}
    n_{j}^{-}> 1, & \hbox{$j=1,...,t_{0}$;} \\
    n_{j}^{-}=0, & \hbox{$j=t_{0}+1,...,t_{0}+t_{1}$;} \\
    n_{j}^{-}=1, & \hbox{$j= t_{0}+t_{1}+1,...,t$.}
  \end{array}
\right.$
\\

\begin{center}

\includegraphics[totalheight=3.5cm]{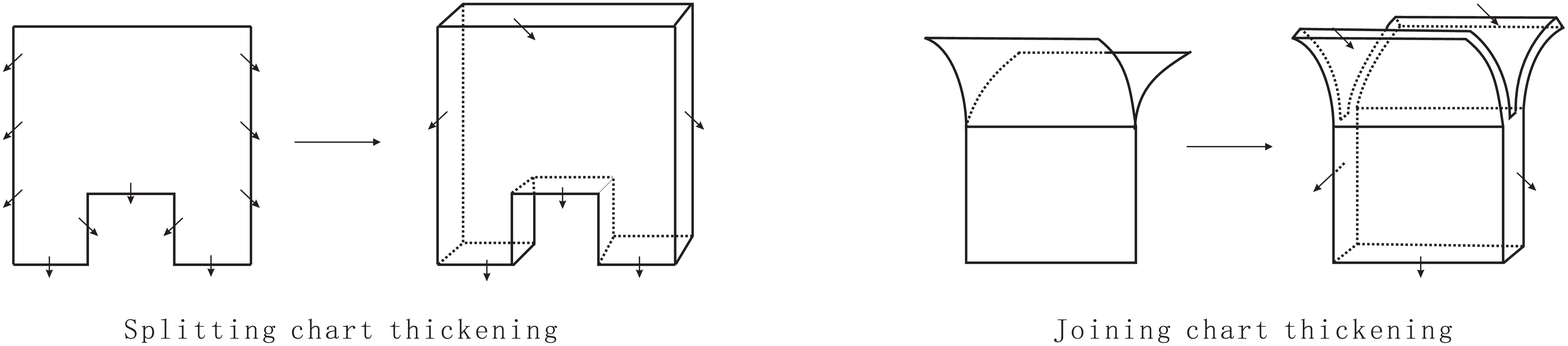}

\begin{center}
Figure 2
\end{center}
\end{center}

\textbf{Definition 2.7} The genus of $T$, $g(T)$ is defined by
$g(T)=max\{\sum_{i=1}^{s_{0}}n_{i}^{+}-s_{0},
\sum_{j=1}^{t_{0}}n_{j}^{-}-t_{0}\}$.

 \vspace*{0.5cm}
\begin{bfseries}
3. The proof of Theorem 1
\end{bfseries}
\vspace*{0.5cm}

The following lemma is Theorem 4.7 in [R]. It is proved by using
Poincare-Hopf formula.

\textbf{Lemma 3.1}
\begin{itshape}
Suppose $\phi_{t}$ is a smooth flow on an odd-dimensional manifold
$M$ which transverses outside to $\partial M^{-}$ and transverses
inside to $\partial M^{+}$, where $\partial M=\partial M^{+} \cup
\partial M^{-}$. Then $\sum I_{p}=\frac{1}{2}(X(\partial
M^{+})-X(\partial M^{-}))$ where the summation is taken over all
singularities in $M$, $I_{p}$ is the index of the singularity $p$
and $X$ denotes the Euler characteristic.
\end{itshape}\\

\textbf{Lemma 3.2}
\begin{itshape}
Let $v$  be a vertex of $L$, which is the Lyapunov graph of an NSF
on a 3-manifold $M$, then $e^{+}-e^{-}-\sum g_{j}^{+}+ \sum
g_{i}^{-}=0$.
\end{itshape}

\textbf{Proof}: To the inverse image (for the map
 $h:M\rightarrow L$) of a small neighborhood of $v$, $\sum I_{p}=0$, $\partial M^{+}=
2e^{+}-2\sum g_{j}^{+}$ and $\partial M^{-}= 2e^{-}-2\sum
g_{i}^{-}$. By Lemma 3.1, we have $0= \frac{1}{2}(2e^{+}-2\sum
g_{j}^{+}-2e^{-}+2\sum g_{i}^{-})$. Therefore
$e^{+}-e^{-}-\sum g_{j}^{+}+ \sum g_{i}^{-}=0$. Q.E.D.\\

\textbf{Lemma 3.3}
\begin{itshape}
For a Lyapunov graph $L$ associated with an NSF on a closed
orientable 3-manifold $M$, if there exists an edge $E\subset L$ such
that its weight is $(n+1), n\geq 0$, then there exists weight 0
edges $E_{1}, ..., E_{n}\subset L$ such that $L$ is connected if we
cut $L$ along $E_{1}, ..., E_{n}$. In particular, $r(L)\geq n$,
where $r(L)$ is the cyclic number of $L$.
\end{itshape}

\textbf{Proof}: Suppose $L(E)$ is  composed of $x\in L$ which
satisfies the condition that there exists an oriented arc starting
at $t(E)$ and ending at $x$, which is a union of a sequence of
oriented edges of nonzero weight.
 Obviously $L(E)$ is a connected
subgraph of $L$. An edge of weight 0 is said to be a \emph{vanishing
weight 0 edge} if its terminal vertex is in $L(E)$. We denote the
set of all vanishing weight 0 edges by $F$.

 For a 3-manifold, if a one
dimensional basic set is a hyperbolic attractor, it must be a closed
orbit attractor. So the terminal edges of $L$ (also $L(E)$) are
weight 1 edges. By Lemma 3.2, for any vertex, $\sum
g_{j}^{+}-e^{+}=\sum g_{i}^{-}-e^{-}$. By this formula and the fact
that a terminal edge of $L$ must be 1, we have the edge number of
$F$,  $\sharp F$, is no less than $n$.

 We obtain a new
graph $L'$ if we cut $L$ along all vanishing weight 0 edges. We
claim that $L'$ is connected. Otherwise, let $L'=L_{1}\sqcup L_{2}$.
Here $L_{1}$ is the connected component containing the subgraph
$L(E)$. The boundaries of $L_{2}$ are exit sets. Let $M_{1}$ and
$M_{2}$ be $h^{-1}(L_{1})$ and $h^{-1}(L_{2})$ respectively. $M_{2}$
admits an NSF which transverses outside to $\partial M_{2}$. Since
we cut $L$ along weight 0 edges, $\partial M_{2}$ is homeomorphic to
a union of $\sharp F$ disjoint 2-spheres. Hence $X(\partial M_{2})=2
\sharp F$. By Lemma 3.1, we have $X(\partial M_{2})=0$. It is a
contradiction, therefore $L'$ is connected.

Since $\sharp F\geq n$, we can choose $E_{1}, ..., E_{n}\subset L$
such that $L$ is connected if we cut $L$ along $E_{1}, ...,
E_{n}$. In particular, the cyclic number of $L$, $r(L)$, is no less than $n$. Q.E.D.\\

\textbf{The proof of necessity in Theorem 1}

\textbf{Proof}: By Lemma 3.3, there exist weight 0 edges $E_{1},
..., E_{n}\subset L$ such that $L$ is connected if we cut $L$ along
$E_{1}, ..., E_{n}$.

   For any given $p_{i}\in int(E_{i}), i=1, ..., n$, $h^{-1}(p_{i})\cong
   S^{2}$. Denote $S^{2}_{i}=h^{-1}(p_{i}), i=1, ..., n$. Obviously they are
   unparallel. $M$ is connected if we cut $M$ along these 2-spheres.
   There exists a neighborhood $W$ of $S^{2}_{1}$ which is
   homeomorphic to $S^{2} \times [0,1]$. Let $V=\overline{M-W}$. $V$
   is connected and
   $\partial V\cong S^{2}\sqcup S^{2}$. Fix two points $p_{1}$ and $p_{2}$ in each of the two
boundary components of V respectively.
   There exists a curve $c$ in $V$ connecting $p_{1}$ with $p_{2}$.
   Choose a cylinder neighborhood  $N(c)$ of $c$. Let $W'=W\cup
   N(c)$, $V'=\overline{M-W'}$. Then $\partial W'=\partial V' \cong
   S^{2}$. Let $M_{1}=V'\cup_{\partial} D^{3}$. $W'\cup_{\partial} D^{3}\cong
   S^{1} \times S^{2}$, so $M=M'\sharp S^{1}\times S^{2}$. Here $D^{3}$ is a 3-ball.
    We can repeat the procedure for each $S_i^2$, to show that $M=M'\sharp nS^{1}\times S^{2}$, where
      $M'$ is a closed orientable 3-manifold. Furthermore $r(L)\geq n$. Q.E.D.\\

   The proof of sufficiency in Theorem 1 uses some constructions which rely on Theorem 2.5.
We first prove some lemmas and propositions.

  \textbf{Lemma 3.4}
  \begin{itshape}
  For any $n\in Z, n\geq 1$, there exists a Smale
  flow $\varphi_{t}$ on $S^{3}$ such that $\varphi_{t}$ satisfies:

(1) Its singularities are composed of  $n$
  singular attractors and $n$ singular repellers;

(2) There exists an $(n+1)T^{2}$ regular level set of the Lyapunov
graph of $\varphi_{t}$.
  \end{itshape}

\textbf{Proof}: It is sufficient to construct a Lyapunov graph $L$
satisfies: (1) There exist $n$ singular attractors vertexes and $n$
singular repellers vertexes in $L$. (2) There are no more vertices
labeled with singularities in $L$. (3) There exists a weight $(n+1)$
edge on $L$. (4) $L$ satisfies the conditions of Theorem 2.5.

    We start with an oriented graph $E$ with only one edge. Suppose the weight of $E$ is $n+1$.
To vertex $i(E)$ ($t(E)$), we add two oriented edges $E_{n}$
($E^{n}$), $F_{n}$ ($F^{n}$) such that $t(E_{n})=i(E),
i(F_{n})=i(E)$ ($i(E^{n})=t(E), t(F^{n})=t(E)$). Here $E_{n}$
($E^{n}$) is a weight $n$ edge and $F_{n}$ ($F^{n}$) is a weight 0
edge.

    Let $A_{1}^{m_{1}}=(a_{ij}^{1}), A_{2}^{m_{2}}=(a_{ij}^{2})$ be
the nonnegative integer irreducible matrices  associated with $i(E),
t(E)$ respectively. Let $B_{l}=(b_{ij}^{l})$, where
$b_{ij}^{l}\equiv a_{ij}^{l}$ $(mod 2)$ and $l=1,2$. Let $k_{l}=\dim
\ker((I-B_{l}):F_{2}^{m_{l}}\rightarrow F_{2}^{m_{l}})$, where
$F_{2}=\mathbb{Z}/2$ and $l=1,2$. By Theorem 2.5, $k_{1}, k_{2}$
satisfy
\begin{equation} \label{eq:1}
\left\{ \begin{aligned}
         k_{1}+1-n \leq 2\leq k_{1}+1\cdots \\
                  k_{1}+1-(n+1)\leq 1\leq k_{1}+1\cdots
                          \end{aligned} \right.
                          \end{equation}

\begin{equation} \label{eq:2}
\left\{ \begin{aligned}
         k_{2}+1-(n+1) \leq k_{1}\leq k_{2}+1\cdots \\
                  k_{2}+1-n\leq 2\leq k_{2}+1\cdots
                          \end{aligned} \right.
                          \end{equation}

We choose $k_{1}=k_{2}=n$. Obviously they satisfy (1) and (2) and
there exist nonnegative integer irreducible matrices $A_{1}, A_{2}$
such that $k_{1}=k_{2}=n$.

Let $n_{1}=n\geq 1$. If $n_{1}> 1$, let $n_{2}=n_{1}-1$ and we add
two oriented edges $E_{n_{2}}$ ($E^{n_{2}}$), $F_{n_{2}}$
($F^{n_{2}}$) such that $t(E_{n_{2}})=i(E_{n_{1}}),
i(F_{n_{2}})=i(E_{n_{1}})$ ($i(E^{n_{2}})=t(E^{n_{1}}),
t(F^{n_{2}})=t(E^{n_{1}})$). Here $E_{n_{2}}$ ($E^{n_{2}}$) is a
weight $n_{2}$ edge and $F_{n_{2}}$ ($F^{n_{2}}$) is a weight 0
edge. Similar to above, we can associate regular matrices to
$i(E_{n_{1}})$ and $t(E^{n_{1}})$. We can repeat the procedure for
$n_{k}$ until $n_{k}=1$. Thus we have constructed a Lyapunov graph
$L$ on $S^{3}$. It is easy to check that it satisfies all conditions
in the first paragraph of the proof. Figure 3 is an example for
$n=2$. Q.E.D.

\begin{center}

\includegraphics[totalheight=3.5cm]{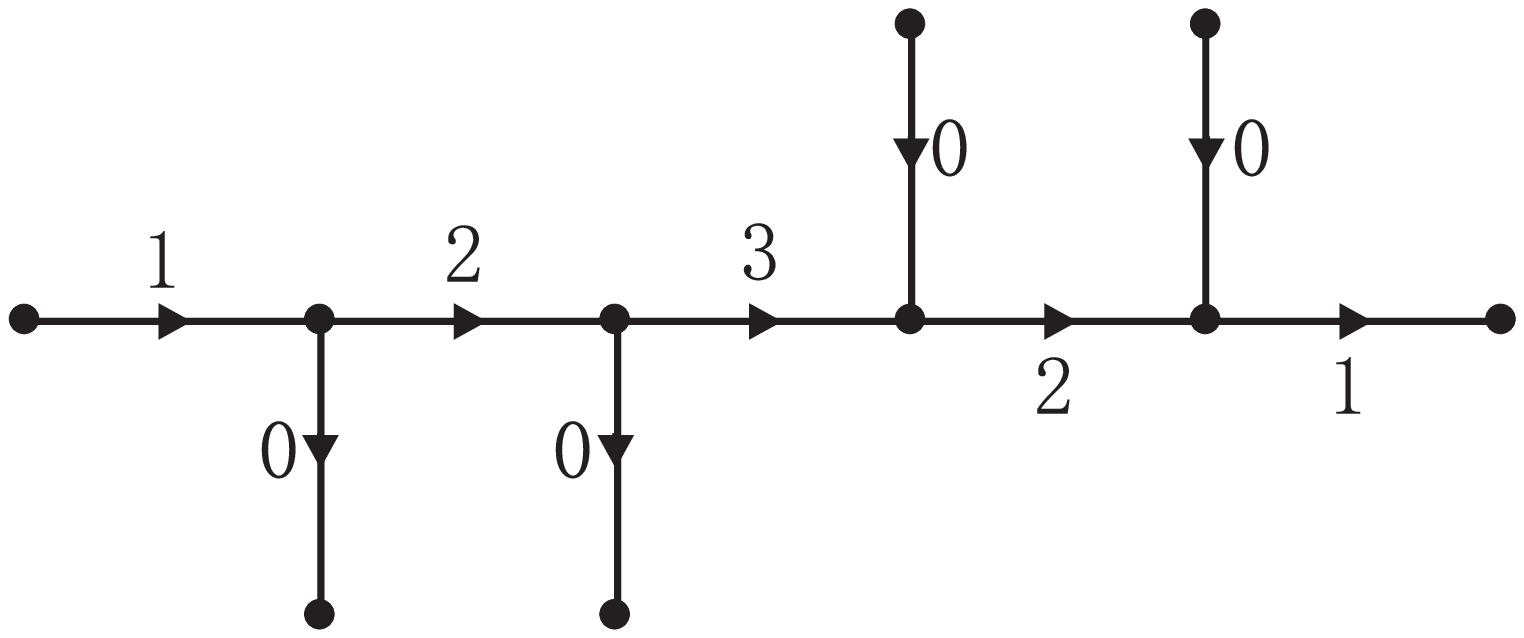}

\begin{center}
Figure 3
\end{center}
\end{center}

\textbf{Proposition 3.5}
\begin{itshape}
For any $n\in Z, n\geq 1$, there exists an NSF $\phi_{t}$ on
$nS^{1}\times S^{2}$ such that there exists an $(n+1)T^{2}$ regular
level set of the Lyapunov graph of $\phi_{t}$.
\end{itshape}

\textbf{Proof}: For any $n\in Z, n\geq 0$, we have constructed a
Lyapunov graph $L$ in the proof of Lemma 3.4. Obviously there is
only one weight $n$ edge on $L$. Suppose the outgoing (incoming)
weight 0 edges in $L$ are $e_{1}, ..., e_{n}$ ($e^{1}, ..., e^{n}$).
Let $x_{i}\in int(e_{i})$ and $y_{i}\in int(e^{i})$, $i=1, ..., n$.

     Cut L along $x_{i}, y_{i}, i=1,...,n$. Denote by $L'$ the connected
component which contains the weight $n$ edge after cutting. Let
     $M=h^{-1}(L')$ and $X_{i}=h^{-1}(x_{i}), Y_{i}=h^{-1}(y_{i}),
     i=1,...,n$. Since $e_{i},e^{i}, i=1,...,n$ are weight 0, $h^{-1}(x_{i}), h^{-1}(y_{i}), i=1,...,n$
     are homeomorphic to $S^{2}$. $\varphi_{t}\mid_{M}$
     transverses outside (inside) to $X_{i}$ ($Y^{i}$), $i=1,...,n$.

     Pasting $X_{i}$ to $Y_{i}$, $i=1,...,n$, we obtain $nS^{1}\times
     S^{2}$. After this surgery,  with some permutations of $\varphi_{t}\mid_{M}$ in a small neighborhood of
     the pasted surfaces, we obtain an NSF $\phi_{t}$ on $nS^{1}\times S^{2}$.
     Since there exists an $(n+1)T^{2}$ regular level set of
      the Lyapunov graph of $\varphi_{t}$ which is disjoint with a small neighborhood of
     the pasted surfaces in $nS^{1}\times S^{2}$, there exists an $(n+1)T^{2}$ regular
    level set of the Lyapunov graph of $\phi_{t}$. Q.E.D.\\

In Figure 4, $G$ is a Lyapunov graph. The two vertexes in $G$
represent a saddle closed orbit and a closed orbit attractor. The
proof of Lemma 3.7 below shows that $G$ is realizable. $L_{1},
L_{2}$ are Lyapunov graphs. $e_{1}, e_{2}$ are weight 1 edges and
$v_{1}, v_{2}$ denote two closed orbit attractors.

 We cut $L_{1}, L_{2}$ along $e_{1},e_{2}$
respectively. Connecting the components of $L_{1}$ and $L_{2}$ which
don't contain $v_{1}$ and $v_{2}$ after cutting with $G$, we obtain
a new Lyapunov graph $L$ denoted by
$L=L_{1}\sharp_{(e_{1},e_{2})}L_{2}$. See Figure 4.

\begin{center}

\includegraphics[totalheight=4cm]{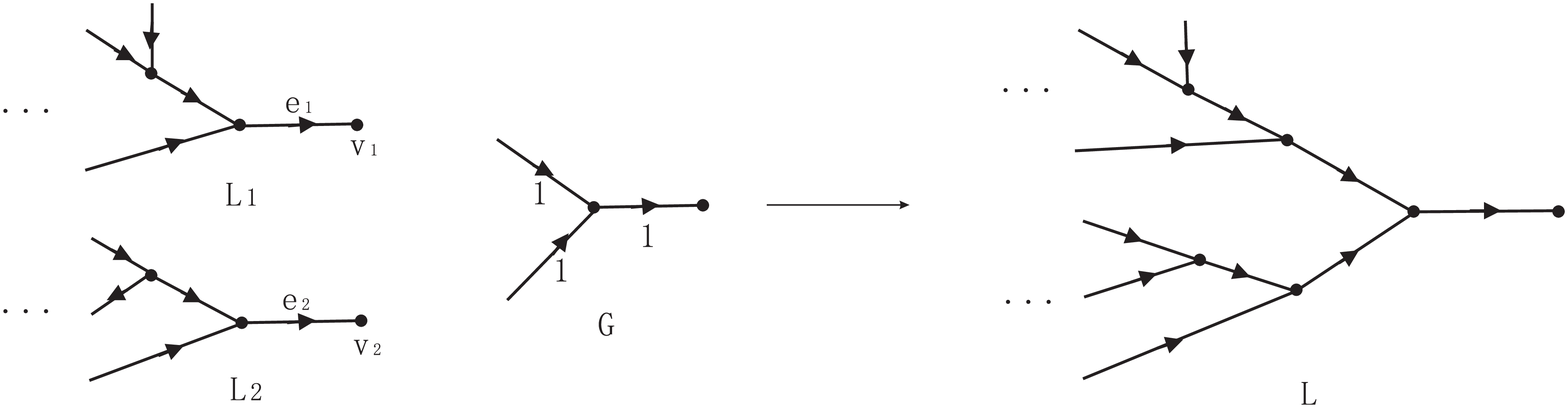}

\begin{center}
Figure 4
\end{center}
\end{center}

\textbf{Lemma 3.6}
\begin{itshape}
Let $M,N$ be two closed orientable 3-manifolds with NSF
$\varphi_{t}^{1}$, $\varphi_{t}^{2}$ respectively and $L_{1}, L_{2}$
are the Lyapunov graphs of $\varphi_{t}^{1}$, $\varphi_{t}^{2}$
respectively. Then there exists an NSF $\phi_{t}$ on $M\sharp N$
such that $L=L_{1}\sharp_{(e_{1},e_{2})} L_{2}$ is the Lyapunov
graph of $\phi_{t}$, where $e_{1}, e_{2}$ are edges of $L_{1},
L_{2}$ respectively which terminate at closed orbit attractor
vertexes.
\end{itshape}

\textbf{Proof}: We will use a surgery method due to Morgan [Mo] in
the proof of this lemma. Before we start the proof, we recall some
basic constructions in [Mo].  Let $W_{1}$ and $W_{2}$ be two solid
tori as shown in Figure 5. There is a flow on $W_{1}$ such that the
center of $W_{1}$ is an attractor and the flow transverse inside to
$\partial W_{1}$. $W_{1}$ with such a flow is called a \emph{round
0-handle}. There is another flow on $W_{2}$ such that the center of
$W_{2}$ is a saddle closed orbit and the flow transverse inside to
$B_{1},B_{2}$ and outside to $\partial M_{2}-\overline{B_{1}\cup
B_{2}}$. $W_{2}$ is one type of the so called \emph{round
1-handles}.

\begin{center}

\includegraphics[totalheight=4cm]{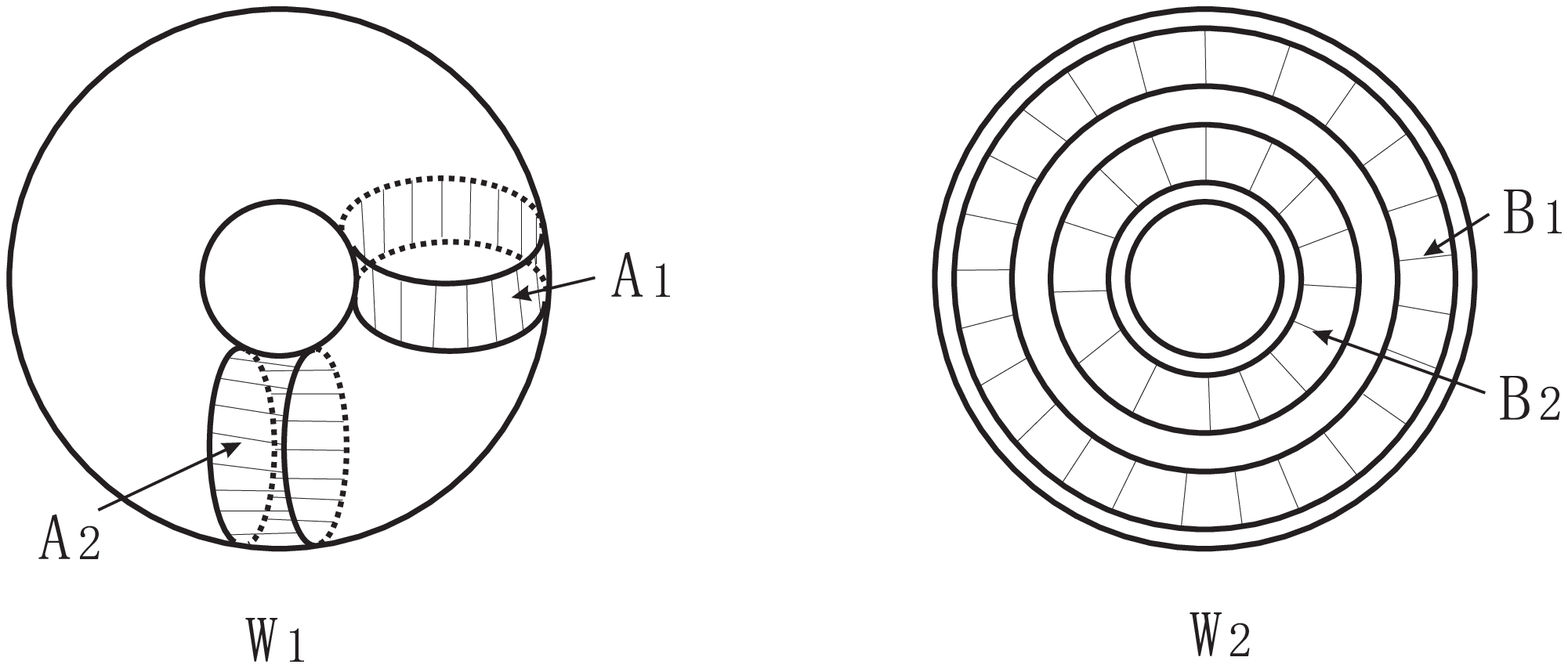}

\begin{center}
Figure 5
\end{center}
\end{center}

$A_{1},A_{2},B_{1}$ and $B_{2}$ are annuluses as shown in Figure 5.
If we attach $A_{1}, A_{2}$ to $B_{1}, B_{2}$ respectively, we
obtain a 3-manifold $W$ with an NSF $\psi_{t}$. $W\cong S^{1}\times
D^{2}\sharp S^{1}\times D^{2}$ and $\psi_{t}$ transverse outside to
$\partial (S^{1}\times D^{2}\sharp S^{1}\times D^{2})$. It is easy
to see that $G$ is a Lyapunov graph of the flow $\psi_{t}$.

  Now, let's turn to the proof of this lemma. Since $\varphi_{t}^{1}$, $\varphi_{t}^{2}$ are NSF on the closed orientable
  3-manifolds $M,N$, there exist closed orbit repellers in $M,N$.
 We cut out two round 2-handles of $M,N$ respectively and call the rest
 parts $M',N'$ respectively. We paste $W$ to $M,N$ along
 their boundaries suitably such that the new  3-manifold is homeomorphic to
 $M\sharp N$. One can check that, $\varphi_{t}^{1}$, $\varphi_{t}^{2}$ and
 $\psi_{t}$ form an NSF $\phi_{t}$ on $M\sharp
 N$ with Lyapunov graph $L=L_{1}\sharp_{(e_{1},e_{2})} L_{2}$.
 Q.E.D.\\

Different from nonsingular Morse-Smale flow (NMSF) and Anosov flow
(For NMSF and Anosov flow, see [Mo] and [Fe]), all closed orientable
3-manifolds admit NSF.

\textbf{Proposition 3.7}
\begin{itshape}
Every closed orientable 3-manifold admits NSF.
\end{itshape}

\textbf{Proof(sketch)}:Proposition 6.1 in [F4] shows that for any
smooth link $K$  in $S^{3}$ there exists an NSF on $S^{3}$ such that
$K$ is the set of attractors of the NSF.

   A famous theorem of Lickorish and Wallace (see, for example, [Rol])
shows that every closed orientable 3-manifold may be obtained by
surgery (cut solid tori and paste solid tori) on a link on $S^{3}$.

   By using the two facts above and some easy combinatorial arguments, one can show
   that every closed orientable 3-manifold admits NSF. Q.E.D.\\

\textbf{The proof of the sufficiency in Theorem 1}

\textbf{Proof}: For any closed orientable 3-manifold $M'$, let
$M=M'\sharp nS^{1}\times S^{2}$ ($n\in \mathbb{N}$). By Proposition
3.7, there exists an NSF $\varphi_{t}$ on $M'$ with a Lyapunov graph
$L_{1}$. By Proposition 3.5 there exists an NSF $\psi_{t}$ on
$nS^{1}\times S^{2}$ with a Lyapunov graph $L_{2}$ such that there
exists an $(n+1)T^{2}$ regular level set and $r(L_{2})\geq n$. By
Lemma 3.6 we have an NSF $\phi_{t}$ on $M=M'\sharp nS^{1}\times
S^{2}$ such that the Lyapunov graph of $\phi_{t}$ is
$L=L_{1}\sharp_{(e_{1},e_{2})} L_{2}$. Furthermore, there exists an
$(n+1)T^{2}$ regular level set in $M$. Since
$r(L_{2})\geq n$, $r(L)\geq n$. Q.E.D.\\

 \vspace*{0.5cm}
\begin{bfseries}
4. The proof of Theorem 2
\end{bfseries}
\vspace*{0.5cm}

\textbf{Lemma 4.1}
\begin{itshape}
For a template $T$, $\sum_{i=1}^{s_{0}} n_{i}^{+} -s_{0}
-s_{1}=\sum_{j=1}^{t_{0}} n_{j}^{-} -t_{0} -t_{1}$
\end{itshape}

\textbf{Proof}: Let $\overline{D}=D\times I$, where $D$ is a disc
and $I=[0,1]$. Pasting  some copies of $\overline{D}$ to
$\overline{T}$ along dividing curves, we have a 3-manifold $M$. The
NSF on $\overline{T}$ can be extended to $M$ in a standard way as an
NSF. By Lemma 3.1, $\sum_{i=1}^{s_{0}}(2-
2n_{i}^{+})+\sum_{s_{0}+1}^{s_{0}+s_{1}}(2-0) $ $=\sum_{j=1}^{t_{0}}
(2-2n_{j}^{-}) +\sum_{t_{0}+1}^{t_{0}+t_{1}}(2-0)$. Therefore
$\sum_{i=1}^{s_{0}} n_{i}^{+} -s_{0} -s_{1}=\sum_{j=1}^{t_{0}}
n_{j}^{-} -t_{0} -t_{1}$. Q.E.D.\\

\textbf{Proposition 4.2}
\begin{itshape}
Suppose $\phi_{t}$ is an NSF on a closed orientable 3-manifold $M$
such that a basic set of $\phi_{t}$ can be modeled by $T$. Then we
can choose a map $h:M\rightarrow L$ with respect to $\phi_{t}$ such
that every $X_{i}$ ($Y_{j}$) is a subset of a regular level set of
$L$. Here $L$ is the Lyapunov graph of $\phi_{t}$.
\end{itshape}

\textbf{Proof}: By the construction of thickened template (See
[GHS]), $\overline{T}$ can be regarded as the suspension of Markov
partition. Then there exists  $N(\overline{T})\supset \overline{T}$
such that $\phi_{t}$ is transverse to $\partial N(\overline{T})$.
Moreover, $N(\overline{T})\cap R(\phi_{t})=\Lambda$. Here $\Lambda$
is the basic set modeled by $T$ and $R(\phi_{t})$ is the chain
recurrent set of $\phi_{t}$. By Theorem 0.3 of [BB],
$N(\overline{T})$ is unique up to topological equivalence.

  [BB] also implies that $N(\overline{T})$ can be obtained from $\overline{T}$
by attaching discs and annuluses along dividing curves. We can
construct a map $h:M\rightarrow L$ with respect to $\phi_{t}$ such
that $\partial N(\overline{T})$ are regular
level sets of $L$. Obviously, $X_{i}(Y_{i})\subset \partial N(\overline{T})$.    Q.E.D.\\

Assume that $\sum_{i=1}^{s_{0}}n_{i}^{+}-s_{0}\geq
\sum_{j=1}^{t_{0}}n_{j}^{-}-t_{0}$. Under this assumption,
$g(T)=\sum_{i=1}^{s_{0}}n_{i}^{+}-s_{0}$. Suppose the conditions of
Proposition 4.2 are true and let $L$ be a Lyapunov graph such that
every $X_{i}$ ($Y_{j}$) is a subset of a regular level set with
respect to $L$. Let $v\in L$ be the vertex modeling $T$ and $e_{1},
..., e_{s}$ be the edges terminating at $v$ whose weight is larger
than 1. The weight of $e_{i}$ ($i=1,...,s$) is denoted by $g_{i}$.

\textbf{Lemma 4.3}
\begin{itshape}
$s\leq s_{0}$ and $\sum_{i=1}^{s}g_{i}-s\geq
\sum_{i=1}^{s_{0}}n_{i}^{+}-s_{0}=g(T)$.
\end{itshape}

\textbf{Proof}: By Proposition 4.2, there exists a map
$P:\{X_{i}\}\rightarrow \{e_{j}\}$, $i=1, ..., s_{0}; j=1, ..., s$.
Since $g_{j}>1$ ($i=1,...,s$), $P$ must be a surjection. Hence
$s\leq s_{0}$. Since $s\leq s_{0}$ and  $\sum_{i=1}^{s}g_{i}\geq
\sum_{i=1}^{s_{0}}n_{i}^{+}$ ($P$ is a surjection),
$\sum_{i=1}^{s}g_{i}-s\geq
\sum_{i=1}^{s_{0}}n_{i}^{+}-s_{0}=g(T)$.   Q.E.D.\\

\textbf{Lemma 4.4}
\begin{itshape}
There exist weight 0 edges $E_{1}, ..., E_{g(T)}\subset L$ such that
$L$ remains connected if we cut $L$ along $E_{1}, ..., E_{g(T)}$. In
particular, $r(L)\geq g(T)$, where $r(L)$ is the cyclic number of
$L$.
\end{itshape}

\textbf{Proof}: The proof is similar to the proof of Lemma 3.3.

Let $L_{1}$ be the set of  $x\in L$ which satisfies the condition
that there exists an arc composed of a sequence of oriented edges of
nonzero weight starting at x and terminating at some $i(E_{i}),
i=1,...,s$. Let $L_{2}= L_{1}\cup e_{1}\cup...\cup e_{s}$. Obviously
$L_{2}$ is a connected subgraph of $L$. A weight 0 edge is said to
be a \emph{vanishing weight 0 edge} if it starts from some vertex in
$L_{1}$. We denote the set of all vanishing weight 0 edges by $F$.

 In 3-manifold, if a
1 dimensional basic set is an attractor, it must be the standard
closed orbit attractor. So those edges of $L_{1}$ (also $L_{2}$)
which are adjacent to a vertex modeling an attractor are weight 1
edges. By Lemma 3.2, For any vertex, $\sum g_{j}^{+}-e^{+}=\sum
g_{i}^{-}-e^{-}$. Due to this formula and the fact that the weight
of the edges above must be 1, we have $\sharp(F)\geq
\sum_{i=1}^{s}g_{i}-s\geq g(T)$ (Lemma 4.3), where $\sharp(F)$ is
the edge number of $F$.

The rest of the proof are the same as the last two paragraphs
 in the proof of Lemma 3.3. Q.E.D.\\

\textbf{The proof of (1) in Theorem 2}

\textbf{Proof}: The proof is the same as the proof of necessity of
Theorem 1 except that we use Lemma 4.4 instead of Lemma 3.3.
Q.E.D.\\

Given $n\in \mathbb{N}$, we cut Lyapunov graph $L$ in the proof of
Lemma 3.4 along the weight $(n+1)$ edge of $L$. $L$ is divided to
two Lyapunov graphs: $L_{1}$, $L_{2}$. $L_{1}$ ($L_{2}$) has weight
$(n+1)$ outside (inside) edge. Let $S^{3}=M_{1}\cup M_{2}$ such that
$M_{i}=h^{-1}(L_{i}), i=1,2$. $M_{1}$ and $M_{2}$ are two genus
$(n+1)$ handlebodies. See [F1] and [Me]. Before we prove part (2) of
Theorem 2, we construct for each $n (n\in \mathbb{N})$ two
handlebodies $H_{1}$ and $H_{2}$ as follows. Cutting small 3-ball
neighborhoods of the singular points in $M_{1}$, $M_{2}$, we obtain
two genus $(n+1)$ handlebodies with $n$ holes, denoted by $H_{1}$
and $H_{2}$ respectively.

\textbf{The proof of (2) in Theorem 2}

\textbf{Proof}: Assume that $\sum_{i=1}^{s_{0}}n_{i}^{+}-s_{0}\geq
\sum_{j=1}^{t_{0}}n_{j}^{-}-t_{0}$. It is easy to show that $X_{i}$
is not a disk for any $i\in \{1,...,s_{0}\}$. We use some annuluses
to connect all $X_{i}, i> s_{0}+s_{1}$ to get a connected surface
$X$. Here the boundaries of the annuluses are attached to some
dividing curves of $\overline{T}$. Obviously the genus of $X$ is 0.
We use an annulus to connect $X$ with $X_{1}$ along a connected
component of $\partial X$ and a connected component of $\partial
X_{1}$. We denote the new surface by $X_{1}'$. The genus of $X_{1}'$
is equivalent to the genus of $X_{1}$, $n_{1}^{+}$. We attach discs
to the boundary of $X_{1}'$ and $X_{i}, 2\leq i\leq s_{0}+s_{1}$.
Then we thicken all attached annuluses and disks to get a new
3-manifold denoted by W. The connected component of $\partial W$
containing $X_1'$ is still denoted by $X_1'$, while the connected
component of $\partial W$  containing $X_i$ ($2\leq i\leq
s_{0}+s_{1}$) is denoted by $X_i'$. Then we extend flow on
$\overline{T}$ to $W$ such that the flow is an NSF. Define the flow
on $W$ by $\varphi_{t}$. For more detail on how to extend the flow,
see Section 3.2 of [Me]. Obviously $\varphi_{t}$ is transverse to
the boundary of $W$.

For any $X_{i}', i> s_{0}$, we attach a round 2-handle to $W$ along
$X_{i}'$. For any $X_{i}', i\leq s_{0}$, if the genus of $X_{i}'$ is
$n+1$, we attach $H_{1}$ to $X_{i}'$. Let $Y_{1}', ..., Y_{t}'$ be
all the outside boundaries of $W$ such that the genus of anyone of
them is larger than 0. We attach $H_{2}$ and round 0-handles to $W$
along $Y_{1}', ..., Y_{t}'$ as above. Therefore, we obtain a
3-manifold $V$ with an NSF $\psi_{t}$. $\partial V = 2g(T)S^{2}$,
$\psi_{t}$ is transverse inside to one copy of $g(T)S^{2}$ and
outside to the other.

Attaching all inside $g(T)S^{2}$ to all outside $g(T)S^{2}$ one by
one, we obtain a 3-manifold $M$ with an NSF $\phi_{t}$. $M=M'\sharp
g(T)S^{1}\times S^{2}$. Here $M'$ is a closed orientable 3-manifold.
Q.E.D.\\

 \vspace*{0.5cm}
\begin{bfseries}
5. A visualization of an NSF on $S^{1}\times S^{2}$
\end{bfseries}
\vspace*{0.5cm}

In the discussion of Theorem 1 and Theorem 2, the existence of NSF
is always proved by some constructions. These constructions are
included in the proofs of the main theorems in [R] and [F1]. [BB]
implies that any NSF on a 3-manifold can be constructed from
thickened templates by some surgeries. In this section, we use
template to construct an NSF on $S^{1}\times S^{2}$ with $2T^{2}$
regular level set.

Denote the template in Figure 6 by $T$. Suppose $\overline{T}$ is
the thickened $T$.

\begin{center}

\includegraphics[totalheight=4.5cm]{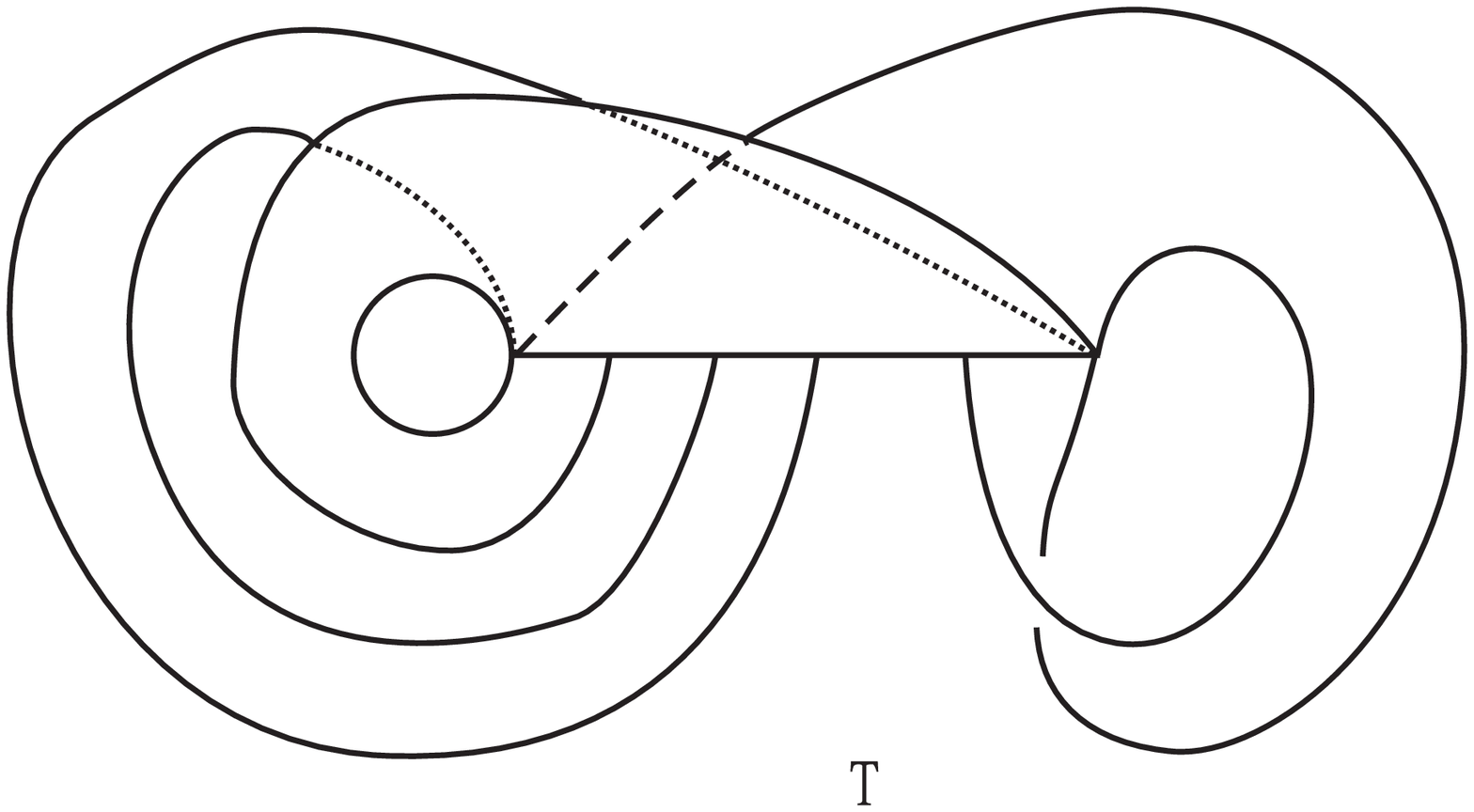}

\begin{center}
Figure 6
\end{center}
\end{center}

The boundary of $\overline{T}$ is composed of entrance sets $X_{1},
X_{2}$ and exit set $Y$. $X_{1}\cap Y=c_{1}$ and $X_{2}\cap
Y=\{c_{2},c_{3},c_{4}\}$ are dividing curves. See Figure 7.

\begin{center}

\includegraphics[totalheight=3.6cm]{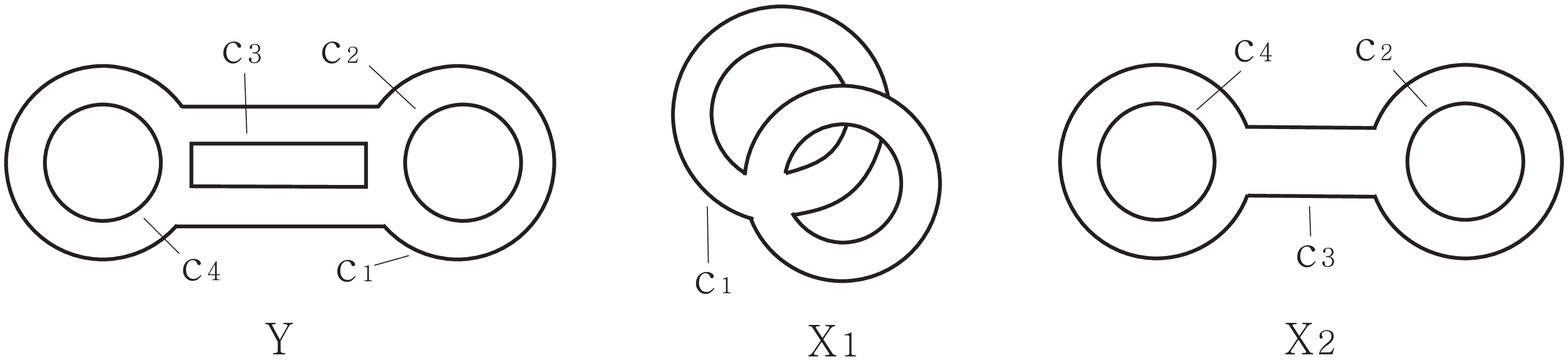}

\begin{center}
Figure 7
\end{center}
\end{center}

We attach three copies of thickened discs $D^{2}\times [0,1]$ to
$\partial \overline{T}$ along $c_{2},c_{3},c_{4}$ (add 2-handles).
We paste thickened punctured torus to $\partial \overline{T}$ along
$c_{1}$. Then we get a new 3-manifold $W$ with an NSF $\varphi_{t}$.
$\partial W=X_{1}'\cup X_{2}'\cup Y'$, $X_{1}'\cong 2T^{2}$,
$X_{2}'\cong S^{2}$ and $Y'\cong T^{2}$. $X_{i}\subset X_{i}',
i=1,2$ and $Y\subset Y'$. $\varphi_{t}$ is transverse inside to
$X_{1}'\cup X_{2}'$ and outside to $Y'$. $W\cong (T^{2}\times
[0,1])\cup_{D_{1}=D_{2}} (S^{1}\times D^{2}- int(D^{3}))$. Here
$D^{i}$ is an i-ball. $D_{1}, D_{2}$ are two discs and $D_{1}\subset
\partial T^{2}\times [0,1]$, $D_{2}\subset
\partial (S^{1}\times D^{2})$. See Figure 8.

\begin{center}

\includegraphics[totalheight=4cm]{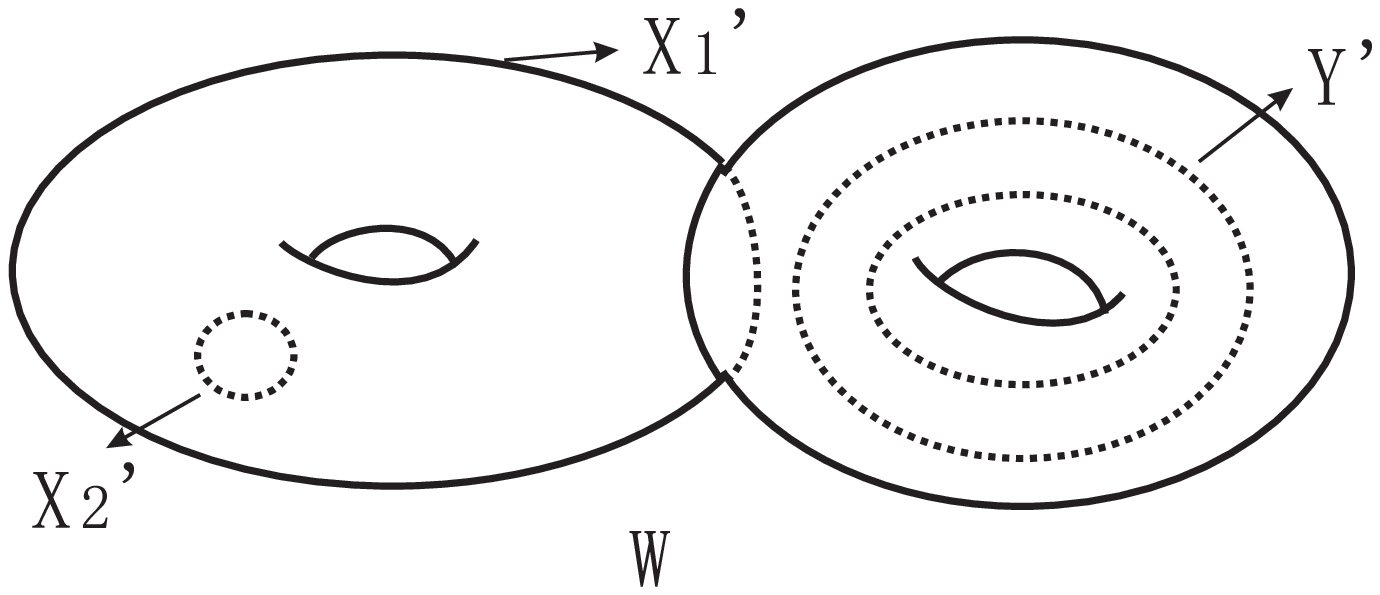}

\begin{center}
Figure 8
\end{center}
\end{center}

The Lyapunov graph of $(W,\varphi_{t})$ is $G$ (Figure 9).  Let
$(V,\psi_{t})$ be the same flow as $(W,\varphi_{-t})$ with an
inverse direction. We add a round 0 handle to $W$ along $Y'$ such
that the new 3-manifold $W'$ is a genus 2 handlebody with a hole.
Similarly, we add a round 2 handle to $V$ along $Y'$ such that the
new 3-manifold $V'$ is a genus 2 handlebody with a hole. Now we
attach $W'$ to $V'$ along their genus 2 surfaces such that the new
3-manifold $M'$ is $S^{3}$ with two holes. $\partial M'=S^{2}\sqcup
S^{2}$. Attaching $M'$ to itself along its two boundary components,
we obtain a 3-manifold $M\cong S^{1}\times S^{2}$ with an NSF
$\phi_{t}$. $L$ (Figure 9) is the Lyapunov graph of $\phi_{t}$.
Obviously there exists a $2T^{2}$ regular level set on $\phi_{t}$.

\begin{center}

\includegraphics[totalheight=3.7cm]{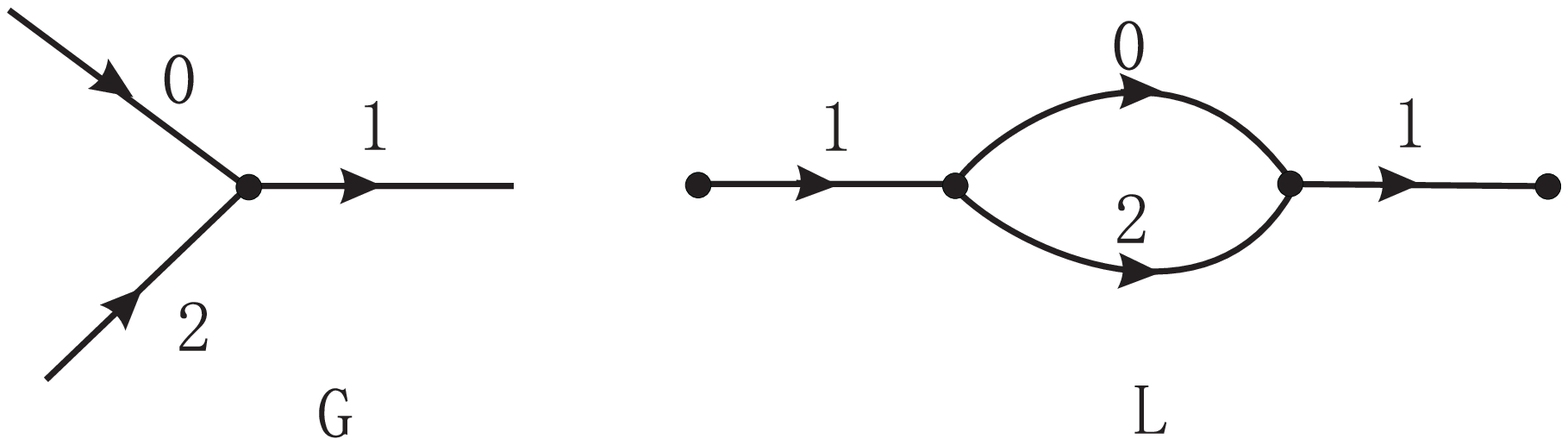}

\begin{center}
Figure 9
\end{center}
\end{center}

Actually, if we don't restrict how to attach $W'$ to $V'$, we can
get any 3-manifold $M\cong M'\sharp S^{1}\times S^{2}$. Here $M'$ is
any closed orientable 3-manifold whose Heegaard genus isn't larger
than 2. Obviously, $M$ admits an NSF and $L$ (Figure 9(2)) is the
Lyapunov graph of the flow.

 \vspace*{0.5cm}
\begin{bfseries}
6. Conclusion
\end{bfseries}
\vspace*{0.5cm}

 Theorem 1 discusses the regular level set of an NSF on a 3-manifold.
The motivation is to generalize J.Franks' work [F1]. The ultimate
goal is to answer,

\textbf{Question 6.1}
\begin{itshape}
How to determine the necessary and sufficient condition when an
abstract Lyapunov graph is associated with an NSF on a 3 manifold?
\end{itshape}\\

Theorem 2 gives some rough description about an NSF on a 3-manifold
with a basic set modeled by a given template $T$. Suppose
$\overline{T}$ is the thickened $T$. Adding 2-handles along all
dividing curves of $\overline{T}$, we obtain a 3-manifold $N(T)$.
[BB] proved that if an NSF on a 3-manifold admits a basic set
modeled by a given template $T$, then there exists a neighborhood of
the basic set which can be obtained by a sequence of standard
surgeries (see [BB]) in $N(T)$. So in some sense understanding the
topological structure of $N(T)$ is not only useful  to our question
but also important to describing NSF on 3-manifolds. Therefore, we
ask,

\textbf{Question 6.2}
\begin{itshape}
What can we say about the topological structure of $N(T)$? For
example, given a template T, is N(T) irreducible?
\end{itshape}

Since the dividing curves of $\overline{T}$ are difficult to
control, describing the topological structure of $N(T)$
systematically is difficult. So discussing whether $N(T)$ is
irreducible  may be more suitable.\\

{\bf Acknowledgements}: The author would like to thank Hao Yin for
his many helpful suggestions and comments. The author also thank M.
Sullivan for telling me the template in Section 5. The author also
would like to thank the referee for carefully reading this paper and
providing helpful comments. The research was supported by the National Natural
Science Foundation of China (grant no. 10926113)\\

{\bf References.} \vskip 4mm

\noindent[BB] B$\acute{e}$guin, F.; Bonatti, C. \emph{Flots de Smale
en dimension 3: pr¨¦sentations finies de voisinages invariants
d'ensembles selles}. (French) [\emph{Smale flows in dimension 3:
finite presentations of invariant neighborhoods of saddle sets}]
Topology 41 (2002), no. 1, 119--162.

\noindent[Bo] Bowen, Rufus. \emph{One-dimensional hyperbolic sets
for flows.} J. Differential Equations 12 (1972), 173--179.

\noindent[BW1] Birman, Joan S.; Williams, R. F. \emph{Knotted
periodic orbits in dynamical systems. I. Lorenz's equations}.
Topology 22 (1983), no. 1, 47--82.

\noindent[BW2] Birman, Joan S.; Williams, R. F. \emph{Knotted
periodic orbits in dynamical system. II. Knot holders for fibered
knots}. Low-dimensional topology (San Francisco, Calif., 1981),
1--60, Contemp. Math., 20, Amer. Math. Soc., Providence, RI, 1983.

\noindent[CR] Cruz, R. N.; de Rezende, K. A. \emph{Cycle rank of
Lyapunov graphs and the genera of manifolds}. Proc. Amer. Math. Soc.
126 (1998), no. 12, 3715--3720.

\noindent[F1] Franks, John. \emph{Nonsingular Smale flows on
$S^{3}$}. Topology 24 (1985), no. 3, 265--282.

\noindent[F2] Franks, John. \emph{Symbolic dynamics in flows on
three-manifolds}. Trans. Amer. Math. Soc. 279 (1983), no. 1,
231--236.

\noindent[F3] Franks, John. \emph{Homology and dynamical systems}.
CBMS49.American Mathematical Society. Providence,Rhode Island,1982.

\noindent[F4] Franks, John.  \emph{Knots,links and symbolic
dynamics.} Ann. of Math. (2) 113 (1981), no. 3, 529--552.

\noindent[Fe] Fenley, S$\acute{e}$rgio R. \emph{Anosov flows in
3-manifolds}. Ann. of Math. (2) 139 (1994), no. 1, 79--115.

\noindent[Fr] Frank, George. \emph{Templates and train tracks}.
Trans. Amer. Math. Soc. 308 (1988), no. 2, 765--784.

\noindent[GHS] Ghrist, Robert W.; Holmes, Philip J.; Sullivan,
Michael C. \emph{Knots and links in three-dimensional flows.}
Lecture Notes in Mathematics, 1654. Springer-Verlag, Berlin, 1997.

\noindent[Mo] Morgan, J. \emph{Nonsingular Morse-Smale flows on
3-dimensional manifolds}. Topology, 18:41-54,1978.

\noindent[Me] V, Meleshuk. \emph{Embedding Templates in Flows}.
Ph.D. dissertation, Northwestern University, 2002.

\noindent[Ok] Oka, Nobuatsu.  \emph{Notes on Lyapunov graphs and
nonsingular Smale flows on three manifolds}. Nagoya Math. J. 117
(1990), 37--61.

\noindent[PS] Pugh, Charles C.; Shub, Michael. \emph{Suspending
subshifts}. Contributions to analysis and geometry (Baltimore, Md.,
1980), pp. 265--275, Johns Hopkins Univ. Press, Baltimore, Md.,
1981.

\noindent[R] de Rezende, Ketty. Smale \emph{flows on the
three-sphere}. Trans. Amer. Math. Soc. 303 (1987), no. 1, 283--310.

\noindent[Ro] Robinson, Clark. \emph{Dynamical systems. Stability,
symbolic dynamics, and chaos}. Second edition. Studies in Advanced
Mathematics. CRC Press, Boca Raton, FL, 1999.

\noindent[Rol] Rolfsen, D. \emph{Knots and Links}. Publish or
Perish, Inc., Berkeley, CA, 1976.


\noindent[Sa] Saito, Masahico. \emph{On closed orbits of Morse-Smale
flows on 3-manifolds}. Bull. London Math. Soc. 23 (1991), no. 5,
482--486.

\noindent[Su] Sullivan, Michael C. \emph{Visually building Smale
flows on $S^{3}$}. Topology Appl. 106 (2000), no. 1, 1--19.

\noindent[Yu] Bin, Yu. \emph{Lorenz like Smale flows on
three-manifolds}. Topology Appl. 156 (2009), no. 15, 2462--2469.

Tongji Department of Mathematics, Tongji University,  Shanghai

200092, P.R.China

binyu1980@gmail.com
\end{document}